\documentclass{article}

\usepackage{amssymb}

\usepackage{amsfonts,amsmath,amsthm}
\usepackage{graphicx,epstopdf}
\usepackage{graphicx}
\usepackage{float}

\usepackage{listings}
\usepackage{color}

\usepackage{float}
\usepackage{amsmath} 
\usepackage{calc}
\newlength{\depthofsumsign}
\usepackage{calc}

\begin{document}



\newcommand{\Eqref}[1]{(\ref{#1})}


\begin{center}
	\large{ \textbf{ {\Large Fast and flexible interpolation via PUM with applications in 
				population dynamics}}}
\end{center}

\begin{center}
	Roberto Cavoretto,  Alessandra De Rossi, Emma Perracchione 
\end{center}

\begin{center}
	Department of Mathematics "G. Peano", University of Turin - Italy
\end{center}
\vskip 0.5cm

\textbf{Abstract.} 	In this paper the Partition of Unity Method (PUM) is efficiently performed using Radial Basis 
Functions (RBFs) as local approximants. In particular, we present a new 
space-partitioning data structure extremely useful in applications 
because of its independence from the problem geometry. Moreover, we study, in 
 the context of wild herbivores in forests, an application of such algorithm.
This investigation shows that the ecosystem of the considered natural park is in a very delicate
situation, for which the animal population could become extinguished.
The determination of the so-called
\emph{sensitivity surfaces}, obtained with the new  fast and flexible interpolation tool,
 indicates some possible preventive measures to the park administrators.

\section{Introduction}
Meshfree methods, first motivated by the fields of geodesy and geophysics, are nowadays 
popular tools for solving PDEs and interpolation problems. 
In this paper,  given scattered samples, the  PUM is performed by blending  RBFs as local approximants.
This problem is  generally considered supposing  to interpolate sets of 
points which are located in suitable  domains, such as the square, \cite{Cavoretto14a,Fasshauer,Wendland05}.
Here, focusing on 2D interpolation, we present a tool extremely useful to work 
in any domain. From this consideration it follows the importance of such a versatile tool in applications.
For instance, it can be employed in population dynamics 
to reconstruct  the   sensitivity surfaces, analyzed in what follows.

As the name of  PUM suggests, an efficient organization of the scattered data among the 
different PU subdomains turns out to be crucial. 
To this aim, in literature, techniques known as \emph{kd-trees}, which are not specifically implemented for PUM, 
have already been designed, \cite{Arya98,Fasshauer}.
In this paper, starting from the results shown in \cite{Cavoretto14a}, we propose a versatile partitioning structure, 
named the \emph{block-based partitioning structure}, 
and a novel related searching procedure extremely suitable for changes in the problem geometry. 
Because of the high flexibility of this splitting-data routine, such method can be employed to approximate 
surfaces defined on irregular domains, as surfaces arising from problems in population dynamics. 
Specifically, an application in the context of wild herbivores in forests is considered, \cite{Fuji,Sabetta,TambVent}. 
We study the case of a natural park located in the  Northern Italy. 
It will be pointed out that it is in a very delicate situation, for which the animal 
population could become extinguished.
The determination of the sensitivity surfaces, via the numerical tool 
presented here, suggests some treatments to prevent the herbivore extinction.

The paper is organized as follows. In Section \ref{interpolazione} we recall theoretical preliminaries on local 
RBF-based PU approximation, while Section
\ref{pum_celle} is devoted to the presentation of the block-based partitioning structure.
In Section \ref{applicazione} an application in population dynamics of such a stretchy algorithm is considered.
Finally Section \ref{conclusioni} deals with conclusions pointing out suggestions to prevent the herbivore extinction.

\section{Partition of unity method via RBF interpolation}
\label{interpolazione}
RBF methods are meshfree and work with data given at scattered node points, \cite{Cavoretto14a,Fasshauer,Wendland05}. 
We consider the case of $N$ distinct points $\{ (x_1,y_1), \ldots , (x_N,y_N)) \} 
\subseteq \Omega$, where $\Omega  \subseteq \mathbb{R}^2$ is an open and bounded domain. 
Let $ \{ f(x_1,y_1), \ldots ,f(x_N,y_N) \}$ be the associated function values, 
then the standard RBF interpolation problem consists in finding an interpolant of the form:
\begin{equation}
P(x,y) = \sum_{i=1}^{N} \lambda_i \varphi(d((x,y),(x_i,y_i))),
\label{rbf1}
\end{equation}
where $d(\cdot,\cdot)$ is the Euclidean distance, 
$\lambda_i \in \mathbb{R}$, for $i= 1, \ldots ,N$ and the function $\varphi$ is a RBF.

The coefficients $\lambda_1, \ldots ,\lambda_N$ are determined
by enforcing the conditions:
\begin{equation*}
 P(x_i,y_i) = f(x_i,y_i), \quad i = 1, \ldots  ,N.
\end{equation*}
 This leads to a
symmetric linear system of equations:
\begin{equation}
A \boldsymbol{\lambda} = \boldsymbol{f},
\label{sys}
\end{equation}
where: 
\begin{equation*}
A_{ij} = \varphi (d((x_i,y_i),(x_j,y_j))), \quad i,j = 1, \ldots N,
\end{equation*}
\begin{equation*}
\boldsymbol{f} = [f(x_1,y_1) , \ldots,  f(x_N,y_N)]^T, 
\end{equation*}  
and 
\begin{equation*}
 \boldsymbol{ \lambda} = [\lambda_1, \ldots, \lambda_N]^T.
\end{equation*}

In order to describe the PUM, let us consider an open and bounded covering of $\Omega$ 
composed by $d$ subdomains or patches $\{\Omega_j\}_{j=1}^{d}$   satisfying  a pointwise
overlap condition. Moreover we also suppose that
$ \forall$ $ (x,y) \in \Omega$, the set $ I(x,y) = \{ j : (x,y) \in \Omega_j \},$ 
is such that $ {\rm card} (I(x,y)) \leq K, $
where the constant $K$ is independent from the number of subdomains.
Associated with the subdomains we choose a 
family of compactly supported, nonnegative and continuous functions $w_j$ subordinate to the subdomain
$\Omega_j$, such that $ \sum_{j = 1}^d w_j( x,y)=1$ on $\Omega$ and 
${\rm supp}(w_j)  \subseteq \Omega_j$. Then the global approximant reads as follows:
\begin{equation} 
P ( x,y)= \sum_{j=1}^{d} P_j( x,y) w_j ( x,y), \quad (x,y) \in \Omega,
\end{equation}
where $P_j$ defines a local RBF interpolant of the form \eqref{rbf1} on 
each subdomain $\Omega_j$ and $w_j: \Omega_j \longrightarrow \mathbb{R}$.
Thus the global interpolant $P$ is obtained by solving $d$ \emph{small} RBF interpolation problems. 
Moreover the global PU approximant inherits the interpolation property of the local interpolants, i.e.:
\begin{align*}
P( x_i,y_i ) = \sum_{j=1}^{d} P_j( x_i,y_i ) w_j (x_i,y_i) = \sum_{j \in I(x_i,y_i)} f( x_i,y_i ) w_j (x_i,y_i) = f( x_i,y_i ).
\end{align*}

\section{The block-based partitioning structure}
\label{pum_celle}
In this section we present our flexible and fast algorithm, built \emph{ad hoc} for bivariate interpolation 
and extremely useful in applications because of its independence from the problem geometry.
We focus on  scattered data points arbitrarily distributed in a (priori unknown) domain $\Omega \subseteq \mathbb{R}^2$. 
Specifically, we approximate the interpolant on the minimal set containing points which
can be automatically detected, i.e. at first we compute the convex hull defined by data sites. 
If the domain is supposed to be known, any generalization is possible and straightforward \cite{Heryudono} and
consequently this strategy is not restrictive in any sense. 

To show our method, used  to quickly store points into the different
subdomains, we first define   the problem geometry. Thus we consider the 
rectangle ${\cal R}$ containing the scattered data: 
\begin{align*}
{\cal R}= [\textrm{min}_x, \textrm{max}_x] \times [\textrm{min}_y, \textrm{max}_y],
\end{align*}
where, as example, with the notation $\textrm{min}_x$ we identify
$\min_i  x_{i}$, $i=1, \ldots, N$.
Moreover, a second auxiliary structure, known 
as \emph{bounding box}, i.e. the square of edge 
$l_{box}= \max  ( \textrm{max}_x, \textrm{max}_y)- \min ( \textrm{min}_x, \textrm{min}_y)$
is also considered.

Then, a uniform grid  of PU subdomain centres ${\cal  C}_{d}= \{(\bar{x}_j, \bar{y}_j )   ,j=1, \ldots, d \}$
is constructed in the convex hull \footnote{This
	set is obtained by generating a grid of $d_{PU}^2$ points on ${\cal R}$, 
	where $d_{PU}= \left( \bigg \lfloor  \frac{\displaystyle 1}{\displaystyle 2} 
	N^{1/2} \bigg \rfloor \right)$. They are then reduced taking only those in $\Omega$.}.
Once the PU structure is built, the problem reduces to solve 
$d$ RBF interpolation problems. 
Consequently, a  partitioning data structure and a related searching procedure 
must be employed to efficiently organize  points among the $d$ subdomains.

In literature, to this scope the  kd-tree partitioning 
structures are commonly and widely used \cite{Arya98,Fasshauer,Wendland05}. 
Here instead, a novel partitioning structure and a related searching procedure, flexible as kd-trees but
with a considerable saving in terms of computational time and cost, is considered.
Specifically, after dividing the bounding box auxiliary structure in 
$q^2$ squares, named blocks, we store points  into such blocks.

The number $q$ of blocks along one side of the bounding box is defined as:
\begin{equation} 
q= \bigg \lceil \frac{\displaystyle l_{box}}{\displaystyle \delta_{PU}} \bigg \rceil,
\label{q}
\end{equation}
where $\delta_{PU}=  l_{box} \sqrt{2} /  d_{PU}$ 
is the radius of the circular patches.
This choice seems to be trivial, but in practice, acting in this way, the search of those
points lying in a subdomain whose centre belongs to the $k$-th block
is limited at the $k$-th neighbourhood, i.e.   the $k$-th block and  its eight neighbouring blocks. 
It leads to a cheap procedure in terms of computational complexity since the searching routine 
is performed in a constant time, independently from the initial number of nodes.


\section{Prevent herbivore extinction via sensitivity surfaces}
\label{applicazione}
In this section it will be pointed out how such a flexible procedure, described
in the previous sections, can be used to determine the sensitivity surfaces, i.e. 
surfaces defined in the space of parameters separating the   stability regions of two different equilibria.

Let $H$, $G$ and $T$ represent respectively the herbivores, grass and trees populations
of the environment in consideration.
The model we consider, \cite{TambVent}, is a classical predator with two prey system,
in which the resources are consumed following a concave response function, usually called
the Beddington-De Angelis function, \cite{De Angelis}. 
Let $\alpha$ and $\beta$ be the inverse of the herbivores maximal consumption of grass  and trees,
respectively.
Letting $r_1$ and $r_2$ denote the grass and trees growth rates, $K_1$ and $K_2$
their respective carrying capacities, $\mu$ the metabolic rate of herbivores,
$c$ and $g$ the half saturation constants,
$e\leq 1$ and $f\leq 1$ the conversion factors of food into new herbivore biomass
and $a$ and $b$ the daily feeding rates due to grass and trees,
respectively, thus the model  reads as follows, \cite{Sabetta}:
\begin{equation}
\begin{aligned}\label{HGTsys}
\dot{H}&=-\mu H +ae\frac{H\ G}{c+H+\alpha\ G}+bf\frac{H\ T}{g+H+\beta\ T+\alpha\ G},\\
\ \\
\dot{G}&=r_1 G\left( 1-\frac{G}{K_1}\right)-a\frac{H\ G}{c+H+\alpha G},\\
\ \\
\dot{T}&=r_2 T\left( 1-\frac{T}{K_2}\right)-b\frac{H\ T}{g+H+\beta\ T+\alpha\ G}.
\end{aligned}
\end{equation}
All parameters are nonnegative; $K_1$, $K_2$, $c$ and $g$ are measured in biomass, 
$e\le 1$, $f\le 1$, $\alpha$ and $\beta$ are pure numbers, $\mu$, $r_1$, $r_2$, $a$ and $b$ are rates.

\begin{figure}[ht]
	\centering
	\includegraphics[height=.22\textheight]{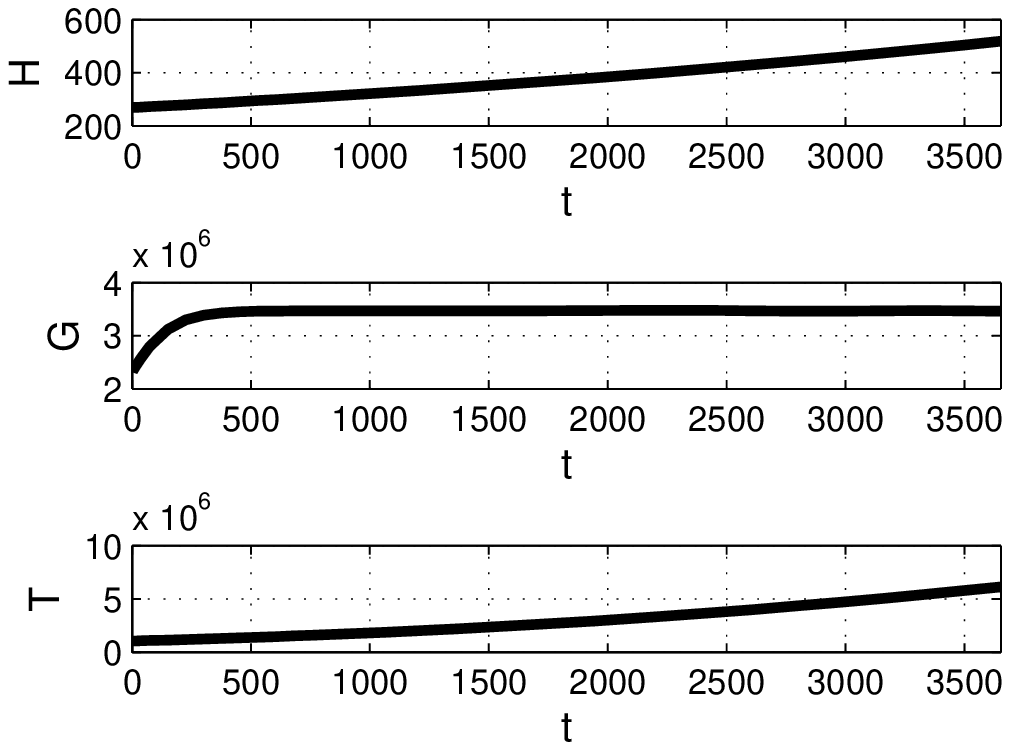}
	\includegraphics[height=.22\textheight]{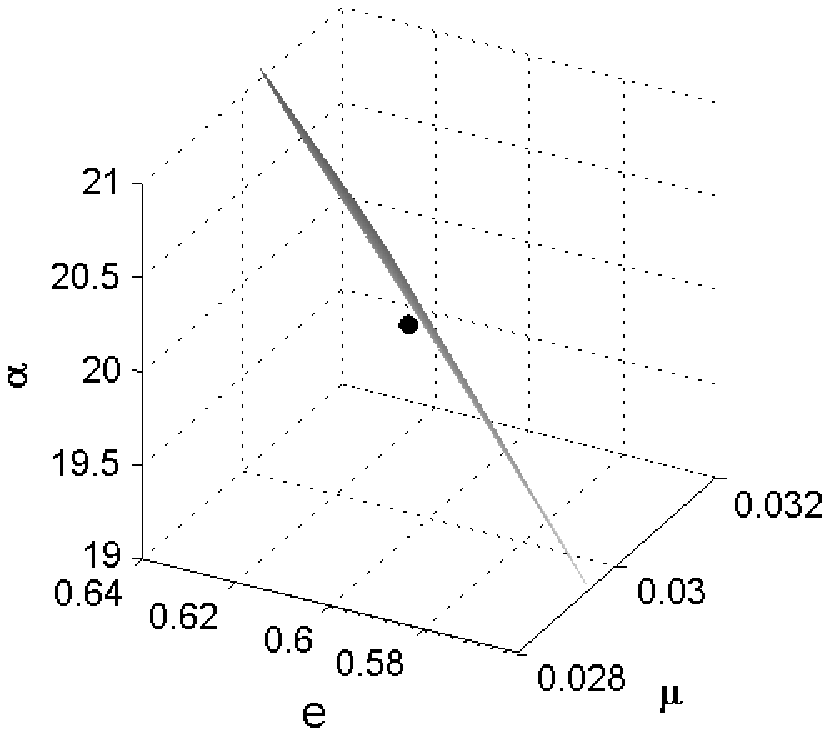} 
	\caption{Left: Dolomiti Bellunesi's system evolution over ten years. 
		Right: the surface separating the stability region of the coexistence equilibrium $E^*$, 		
		bottom, from the one of the
		herbivore-free equilibrium $E_1$, top right. The dot represents the actual situation.
	}
	\label{bellunesi_10years}
\end{figure}

For the study of critical points and stability of System \eqref{HGTsys} refer to \cite{TambVent}.
The equilibria, which play a role in this investigation, are 
the herbivore-free equilibrium $E_1=(0,K_1,K_2)$ and the coexistence equilibrium point $E^*=(H^*,G^*,T^*)$.
The latter can be assessed only via numerical simulations.
We test the model considering the natural park of the Dolomiti Bellunesi located in the  Northern Italy. 
To estimate the parameters we refer to \cite{Fuji,TambVent}. 
Basing our considerations on data related to the number of herbivores, the extension of the park and 
following  tables providing the estimation of the annual net primary production of several
environments  shown in \cite{Fuji}, we can fix $\mu=0.03$,  $r_1=0.01$, $r_2=0.0006$, $\alpha=0.05^{-1}$, $\beta=8$,
$e=0.605$, $f=0.001$, $K_1={3469640.64}$, $K_2={15695993.39}$, $c={101862.16}$ and $g={1001229580.18}$.
Moreover, also the initial conditions $H(0)=268.750$, $G(0)=2313093.76$ and $T(0)=1046399.56$ are known.
With such parameters and initial conditions the model behaves as shown in Figure \ref{bellunesi_10years} (left).
But we noticed by mean of numerical simulations that the herbivore population level appears
to be very sensitive to small perturbations of several parameters and under the  high risk of extinction.
Figure \ref{bellunesi_10years} (right), in the three-dimensional parameter space $(\mu,e,\alpha)$, represents
the sensitivity surface separating the stability region of the coexistence equilibrium $E^*$ from the one 
of the herbivore-free equilibrium $E_1$. 

To reconstruct such surface we consider as RBF 
the compactly supported Wendland's $C^2$ function $
\varphi(r) =  \left( 1-\varepsilon r \right)_+^{4} \left(4\varepsilon r +1\right), $
where $\varepsilon  \in  \mathbb{R}^{+}$ is the so-called
shape parameter and $(\cdot)_+ $ denotes the truncated power function.  
Moreover, since the sensitivity surface is  defined on a convex domain, it has been reconstructed
via the versatile and flexible partitioning 
structure described in Section \ref{pum_celle}. The points describing  
the surface have been found with a bisection-like routine similar to the one 
outlined in \cite{C-D-P-V}. 

We have already pointed out the versatility of the partitioning structure adopted in this paper.
Concerning the efficiency, the CPU times, required  to interpolate the 
$N=107$ \emph{sensitivity points},  with the block-based data structure and, for comparison, with the use of  
kd-trees\footnote{We use the only full \textsc{Matlab} implemented package for kd-trees, written by P. 
	Vemulapalli, available at \textsc{Matlab} central file exchange \emph{http://www.mathworks.com/matlabcentral/fileexchange/}} are
$t_{bb}=0.31s$ and $t_{kd}=5.02s$, respectively.

\section{Final remarks}
\label{conclusioni}
We presented a new fast and flexible interpolation tool and an application to herbivore population dynamics. From Figure \ref{bellunesi_10years} (right) it is easy to verify that the actual value of 
the herbivore population, represented by the dot,
is very close to the separatrix surface. Thus  changes in the environmental
conditions that might lead even to relatively small parameter perturbations may 
push it into the region where the herbivores would be doomed.

Moreover, the representation of this sensitivity surface enables us to 
suggest some treatments in order to possibly prevent the herbivore extinction.
Mathematically speaking, the strategy consists in moving away from the sensitivity surface
the current value of the parameters $\mu$, $e$ and $\alpha$.
For instance, recalling that $\mu$ is the mortality rate, building 
safety niches during  winter surely leads to a decrease of $\mu$.
Furthermore, by planting grass more nutritious, the nutrient assimilation coming from grass, i.e.
$\alpha^{-1}$, increases.

\end{document}